\documentclass{amsart}
\usepackage{amssymb,amsmath,epsfig,mathrsfs, enumerate}
\usepackage[normalem]{ulem}
\usepackage{graphicx}
\usepackage[normalem]{ulem}
\usepackage[utf8]{inputenc}  
\usepackage[T1]{fontenc}     
\usepackage[margin=3cm]{geometry}
\usepackage{hyperref}
 \usepackage[pagewise]{lineno}
\hypersetup{
 colorlinks   = true,
 urlcolor     = blue,
 citecolor   = red ,
 bookmarksopen=true
}
\newtheorem{theorem}{Theorem}[section]
\theoremstyle{plain}

\newtheorem{definition}[theorem]{Definition}

\newtheorem{lemma}[theorem]{Lemma}

\newtheorem{remark}[theorem]{Remark}

\numberwithin{equation}{section}
\title[Multiplicity results in for a class ]{Multiplicity Results for Fully Nonlinear Elliptic Equations with Natural Gradient Growth}

\author{Mohan Mallick}
\address[Mohan Mallick]{VNIT Nagpur, India-440010} \email{mohan.math09@gmail.com, mohanmallick@mth.vnit.ac.in}

\author{Ram Baran Verma}
\address[Ram Baran Verma]{SRM University Amaravati, Andhra Pradesh-522502, India} \email{rambaran.v@srmap.edu.in,rambv88@gmail.com}
\date{}
\begin{document}
\maketitle
\begin{abstract}
\noindent In this paper, we prove a theorem concerning the existence of three solutions for the following boundary value problem:
\begin{equation*}
-\mathcal{M}_{\lambda,\Lambda}^+(D^2u)-\Gamma|Du|^2=f(u)~~~\text{in}\ \Omega,
u=0~~~\text{on}\ \partial\Omega,
\end{equation*}
where $f:[0,\infty]\to[0,\infty]$ is a $C^{\alpha}$ function and $\Omega$ denotes a bounded, smooth domain in $\mathbb{R}^N$. By constructing two ordered pairs of sub and supersolutions for a specific class of $f$ exhibiting sublinear growth, we further establish the existence of three positive solutions to the aforementioned boundary value problem.
\end{abstract}\noindent{\bf Mathematics Subject Classification (2010):} { Primary {35J25, 35J60}; Secondary 35D40, 35A01}.
\section{Introduction}
In this article, we investigate the existence of multiple solutions, specifically three positive solutions, to the following boundary value problem:
\begin{equation}\label{problem}
\left\{
\begin{aligned}
-\mathcal{M}_{\lambda,\Lambda}^+(D^2u) - \Gamma |Du|^2 &= f(u) \quad \text{in} \quad \Omega,\\
u &= 0 \quad \text{on} \quad \partial\Omega,
\end{aligned}
\right.
\end{equation}
where $\mathcal{M}_{\lambda,\Lambda}^+$ denotes the Pucci extremal operator, $\Omega$ represents a bounded domain with a smooth boundary in $\mathbb{R}^N$ for $N \geq 1$, and $f$ satisfies the following condition:

\textbf{(C0)} $f \in C_{loc}^{\alpha}[0,\infty)$ for some $\alpha > 0$, and there exists $\gamma > 0$ such that $f(r) - f(s) > -\gamma(r - s)$ for any real $r, s > 0$ with $r > s$.

In seeking to contextualize our research within the existing body of work, it is instructive to review key contributions made by researchers in this field. We begin by observing that the second-order operator featured in Problem \ref{problem} is non-variational, and the corresponding equation behaves properly when the function \( f \) is non-increasing in its arguments. In mathematical terms, a proper operator is analogous to a coercive operator within a variational framework.

The exploration of solutions to equations characterized by natural gradient growth commenced with the seminal work of L. Boccardo et al.\cite{boccardo1982existence} in 1982. Subsequently, numerous scholars have engaged with this topic, addressing various complexities associated with natural gradient growth, as documented in works such as \cite{boccardo1983existence, bensoussan1988non, maderna1992quasilinear, ferone2000nonlinear, boccardo1992estimate}. For insights into fourth-order equations, readers are directed to \cite{dwivedi2023biharmonic}.

These studies predominantly focus on operators with variational and coercive properties. Furthermore, the comparison principle, a critical aspect of these studies, has been extensively examined in literature such as \cite{barles1995uniqueness, barles2006uniqueness, barles1999remarks}. Notably, the initial findings on the existence of solutions employing noncoercive operators were presented in \cite{hamid2008connection}, marking a significant expansion in the understanding of such mathematical structures.

The existence of multiple solutions to coercive equations with variational structure has been investigated by many authors \cite{jeanjean2013existence,de2018existence,abdellaoui2017multiplicity,arcoya2015continuum,de2017multiplicity,de2019priori,chaouai2023priori}. The first attempt to establish the existence of multiple solutions to noncoercive equations was made in \cite{jeanjean2013existence}, with subsequent multiplicity results for noncoercive operators in \cite{de2017multiplicity}. The approach in this paper is adaptable to equations with a non-variational structure. Multiplicity results for quasilinear equations with natural growth have been established in \cite{chaouai2023priori,de2020existence,rendon2022multiplicity}.

Fully nonlinear elliptic equations of the nonvariational form have also attracted continuous attention. The existence of solutions to equations with natural growth in gradient was first addressed in \cite{MHP,ishii1990viscosity}. The existence and regularity properties of solutions to equations with non-smooth coefficients have been studied in \cite{koike2004maximum,koike2009weak}. In \cite{sirakov2010}, the author extensively studied problems involving the Pucci operator, establishing existence and global $C^{\alpha}$ regularity of solutions. The existence of positive solutions to similar equations with singular nonlinearity has been proved in \cite{tyagi2019positive,birindelli2019dirichlet}. Recently, \cite{Sirakovm6} proved the existence of at least two positive solutions to equations of the type (\ref{problem}) with superlinear growth in gradient and linear-type $f$. For further references, see \cite{PA3,JT}.

In this setting, our problem is non-proper as we assume $f$ to be non-negative and non-decreasing. Multiplicity of solutions is naturally exhibited in such equations \cite{sirakov2010}. Our results extend those of \cite{Sirakovm6} by obtaining three solutions. While our equations have constant coefficients, unlike the more general coefficients in \cite{Sirakovm6}, our existence results do not rely on the assumption of a solution in the borderline case ($\mu=0$). We also construct sub- and super-solutions explicitly.

This paper focuses on two aspects: the development of the theory related to the existence of three positive solutions and an example fitting our setting. The idea of establishing three positive solutions first appeared in the work of Amann \cite{aman} and Shivaji \cite{shivaji1987remark} for the Laplace operator.

The paper is structured into five sections. Section 2 introduces definitions and preliminary results related to the Pucci extremal operator. Section 3 is dedicated to proving the theorem that establishes three solutions for Problem (\ref{problem}). In Section 4, we state and demonstrate the existence and multiplicity results for a problem of the type (\ref{problem}), characterized by sublinear growth at infinity. The final section presents an example of a sublinear function $f$ that allows Problem (\ref{problem}) to exhibit at least three positive solutions.

\section{Preliminaries}

In this section, we present the essential definitions and results that will be used throughout the article. For specified constants where $0 < \lambda < \Lambda$, Pucci's extremal operators are defined as follows:
\begin{equation}\label{pucci}
\mathcal{M}^{\pm}_{\lambda,\Lambda}(X)=\Lambda\sum_{\pm e_{i}>0}e_{i}+\lambda\sum_{\pm e_{i}<0}e_{i},
\end{equation}
where $X$ is a symmetric $N\times N$ matrix, and $e_i$ are its eigenvalues. The operator under consideration has a non-divergent structure, which leads to a different notion of weak solutions, known as "Viscosity solutions."

For a more general problem than \eqref{problem}, consider:
\begin{equation}\label{P1}
-\mathcal{M}_{\lambda,\Lambda}^{+}(D^2 u)-\Gamma|Du|^{2}+\gamma u=h(x).
\end{equation}

\begin{definition}[see \cite{MHP}]\label{viscosity}
A function $u \in C(\bar{\Omega})$ qualifies as a \emph{viscosity subsolution} (resp., \emph{supersolution}) of \eqref{P1} if, for any test function $\phi \in C^{2}(\Omega)$, $u-\phi$ attains a local maximum (resp., minimum) at a point $x_0 \in \Omega$. At such a point $x_0$, the function $u$ must satisfy the inequality:
\begin{equation*}
-\mathcal{M}_{\lambda,\Lambda}^+ (D^2\phi(x_0)) - \Gamma |D\phi(x_0)|^2 + \gamma u(x_0) \leq h(x_0) \quad (\text{resp.,} \geq h(x_0)).
\end{equation*}
A function $u$ is recognized as a solution if it simultaneously serves as both a subsolution and a supersolution. Moreover, $u$ is defined as a \emph{strict subsolution} (resp., \emph{supersolution}) of \eqref{P1} if it fulfills the subsolution (resp., supersolution) criteria for:
\begin{equation}
-\mathcal{M}_{\lambda,\Lambda}^+ (D^2u) - \Gamma|Du|^2 + \gamma u = h(x) + g(x),
\end{equation}
where $g \in C(\bar{\Omega})$ exhibits the property that $g(x) < 0$ (resp., $g(x) > 0$) throughout $\bar{\Omega}$.
\end{definition}

As discussed in the introduction, the comparison principle for equations with natural growth in gradient has been established by Ishii and Lions \cite{ishii1990viscosity}. However, we adopt the following theorem concerning the comparison principle from \cite{sirakov2010}:

\begin{theorem}\label{comp}[Proposition 3.1\cite{sirakov2010}]
Let $u,v\in C(\overline{\Omega})$ satisfy:
\begin{equation*}
\left\{
\begin{aligned}
-\mathcal{M}_{\lambda,\Lambda}^+(D^2u)-\Gamma|Du|^2+\gamma u&\leq h(x) &&\text{in } \Omega,\\
-\mathcal{M}_{\lambda,\Lambda}^+(D^2v)-\Gamma|Dv|^2+\gamma v&\geq h(x) &&\text{in } \Omega,\\
u&\leq v &&\text{on } \partial\Omega,
\end{aligned}
\right.
\end{equation*}
then $u\leq v$ in $\overline{\Omega}$.
\end{theorem}

To demonstrate the multiplicity of solutions, we require an estimate on the gradient of the solutions to this class of equations. The following result provides an a priori $C^{1,\alpha}$-estimate for the considered class of equations, where the structure conditions $(SC)^{\mu}$ and $H_{\theta}$, as mentioned in \cite{nornberg2019}, are satisfied:

\begin{theorem}\label{c1alpha}[Theorem 1.1\cite{nornberg2019}]
Let $u$ satisfy:
\begin{equation*}
\left\{
\begin{aligned}
-\mathcal{M}_{\lambda,\Lambda}^+(D^2u)-\Gamma|Du|^2+\gamma u&= f &&\text{in } \Omega,\\
u&=0 &&\text{on } \partial\Omega,
\end{aligned}
\right.
\end{equation*}
with $\|u\|_{L^{\infty}(\Omega)}+\|f\|_{L^{p}(\Omega)}\leq C_1$ for $p>n$. Then there exists $\alpha\in(0,1)$ depending on $N, p, \lambda, \Lambda,$ such that $u\in C^{1,\alpha}(\overline{\Omega})$ and satisfies the estimate:
\begin{equation}
\|u\|_{C^{1,\alpha}(\overline{\Omega})}\leq C\Big(\|u\|_{L^{\infty}(\Omega)}+\|f\|_{L^{\infty}(\Omega)}\Big),
\end{equation}
where $C$ depends on $n, \lambda, \Lambda, \Gamma, \gamma, C_1$, and the diameter and smoothness of $\Omega$.
\end{theorem}

In the above theorem, we have used the condition $\|u\|_{L^{\infty}(\Omega)}+\|f\|_{L^{\infty}(\Omega)}\leq C_1$. This will be achieved with the help of the comparison principle combined with the following proposition:

\begin{lemma}[Lemma 3.1 \cite{sirakov2010}]
For any positive constants $\gamma, \Gamma,$ and $k,$ there exist viscosity solutions $u_1$ and $u_2$ satisfying:
\begin{equation}
\label{lll}
\left\{
\begin{aligned}
-\mathcal{M}_{\lambda,\Lambda}^+(D^2u_1)-\Gamma|Du_1|^2+\gamma u_1&\geq k &&\text{in } \Omega,\\
-\mathcal{M}_{\lambda,\Lambda}^-(D^2u_2)-\Gamma|Du_2|^2+\gamma u_2&\leq -k &&\text{in } \Omega,\\
u_2\leq 0 &\leq u_1 &&\text{in } \Omega,\\
u_2= u_1&=0 &&\text{on } \partial\Omega,
\end{aligned}
\right.
\end{equation}
and $\|u_{i}\|_{L^{\infty}(\Omega)}\leq \Big(\frac{\lambda}{\Gamma}\Big)\Big(e^{\frac{\Gamma k}{\lambda\gamma}}-1\Big)$ for $i=1,2.$
\end{lemma}

\subsection{Compactness of the Solution Operator}

From the hypothesis (C0), we know that $f(s)+\gamma s$ is increasing in $s$. For a given $u\in C_0(\bar{\Omega})$, consider the unique solution $w\in C_{0}(\bar{\Omega})$ to the following boundary value problem: 
\begin{equation}\label{existence1}
\left\{
\begin{aligned}
-\mathcal{M}_{\lambda,\Lambda}^+(D^2w)-\Gamma|Dw|^2+\gamma w&=f(u)+\gamma u &&\text{in } \Omega,\\
w&=0 &&\text{on } \partial\Omega.
\end{aligned}
\right.
\end{equation}
For the existence and uniqueness of $w$, see Theorem 1 and Corollary 3.1 in \cite{sirakov2010}. Now we define the operator $S_{f}:C_{0}(\bar{\Omega})\rightarrow C_{0}(\bar{\Omega})$ as $S_{f}(u)$ being the unique solution to \eqref{existence1} for a given $u\in C_{0}(\bar{\Omega})$.

\textbf{$S_{f}$ is Compact:} Let $\{u_n\}$ be a bounded sequence in $C_{0}(\bar{\Omega})$, meaning there exists an $M_1$ independent of $n$ such that $\|f(u_n)+\gamma u_n\|_{L^{\infty}(\Omega)}\leq M_1$. Let $w_n=S_{f}(u_n)$. Observe that $w_n$ satisfies:
\begin{equation}
\label{existence}
\left\{
\begin{aligned}
-\mathcal{M}_{\lambda,\Lambda}^+(D^2w_n)-\Gamma|Dw_n|^2+\gamma w_n&\leq M_1 &&\text{in } \Omega,\\
-\mathcal{M}_{\lambda,\Lambda}^-(D^2w_n)-\Gamma|Dw_n|^2+\gamma w_n&\geq -M_1 &&\text{in } \Omega,\\
w_n&=0 &&\text{on } \partial\Omega.
\end{aligned}
\right.
\end{equation}
Let $u_1$ and $u_2$ be the solutions of \eqref{lll} with $k=M_1$. Then, by the comparison principle, we have $u_2\leq w_n\leq u_1$, and thus:
\begin{equation}\label{lll2}
\|w_n\|_{L^{\infty}(\Omega)}\leq \Big(\frac{\lambda}{\Gamma}\Big)\Big(e^{\frac{\Gamma M_{1}}{\lambda\gamma}}-1\Big).
\end{equation}
Hence, we have:
\begin{equation}
\|w_n\|_{L^{\infty}(\Omega)}+\|f(u_n)+\gamma u_n\|_{L^{\infty}(\Omega)}\leq \Big(\frac{\lambda}{\Gamma}\Big)\Big(e^{\frac{\Gamma k}{\lambda\gamma}}-1\Big)+M_1.
\end{equation}
Using Theorem \ref{c1alpha}, we obtain:
\[\|w_n\|_{C^{1,\alpha}(\Omega)}\leq C\Big(\|w_n\|_{L^{\infty}(\Omega)}+\|\gamma u_n +f(u_n)\|_{L^{\infty}(\Omega)}\Big)\leq C_1,\]
where the compact embedding of $C^{1,\alpha}(\Omega)$ into $C^{1}(\overline{\Omega})$ implies the existence of a subsequence of $\{w_n\}$ that converges in $C^{1}(\overline{\Omega})$ and, therefore, in $C_0(\overline{\Omega})$.

\textbf{Continuity of $S_f$:}
Suppose $\{u_n\}\subset C_0(\overline{\Omega})$ such that $u_n\rightarrow u$ in $C_0(\overline{\Omega})$. From the above calculation, we find that $S_f(u_n)$ converges to some $U$ in $C^{1}(\overline{\Omega})$. Our claim is to show that $U=S_{f}(u)$. In view of the comparison principle (Uniqueness), it suffices to show that $U$ satisfies:
\begin{equation}\label{conv}
\left\{
\begin{aligned}
-\mathcal{M}_{\lambda,\Lambda}^+(D^2U)-\Gamma|DU|^2+\gamma U&=f(u)+\gamma u &&\text{in } \Omega,\\
U&=0 &&\text{on } \partial\Omega.
\end{aligned}
\right.
\end{equation}
Taking $\phi\in C_{\text{loc}}^{2}(B)$, where $B\subset\subset\Omega$, consider:
\[g(x)=-\mathcal{M}_{\lambda,\Lambda}^+(D^2\phi)-\Gamma|D\phi|^2+\gamma U-f(u)-\gamma u,\]
and
\[g_{n}(x)=-\mathcal{M}_{\lambda,\Lambda}^+(D^2\phi)-\Gamma|D\phi|^2+\gamma S_{f}(u_{n})-f(u_{n})-\gamma u_{n}.\]
Then, we have:
\[g_{n}(x)-g(x)=\gamma S_{f}(u_{n})-f(u_{n})-\gamma u_{n}-\gamma U+f(u)+\gamma u.\]
Since $u_n$ converges to $u$ uniformly and $f$ is locally H\"{o}lder continuous, we get:
\[\|g_{n}-g\|_{L^{\infty}(\Omega)}\rightarrow0.\]
By the stability result for viscosity solutions (Proposition 2.3 in \cite{nornberg2019}), we find that $U$ solves \eqref{conv}. Therefore, by the uniqueness of the solution to \eqref{conv}, we conclude that $U=S_{f}(u)$. This result holds for all subsequences of $\{u_n\}$, hence $S_{f}(u_n)$ converges to $U$ in $C^{1}(\overline{\Omega})$.

\begin{remark}\label{rem11}
The above theorem shows that if $u_{n}$ converges uniformly to $u$, then $S_{f}(u_{n})$ converges to $S_{f}(u)$ in $C^{1}(\bar{\Omega})$.
\end{remark}

Before proceeding further, we introduce a new function space $C_{e}(\overline{\Omega})$, where $e$ is the solution to the following boundary value problem:

\begin{theorem}[Theorem 17.18\cite{GT}]\label{Pr2.6}
Let \( e \) in \( C^2(\Omega) \cap C(\bar{\Omega}) \) be the unique solution to the problem:
\begin{equation}\label{ee}
\left\{
\begin{aligned}
-\mathcal{M}_{\lambda,\Lambda}^+ (D^2e) &= 1 &&\text{in } \Omega,\\
e &= 0 &&\text{on } \partial \Omega.
\end{aligned}
\right.
\end{equation}
\end{theorem}

It is established that \( e \geq 0 \). Applying the strong maximum principle and H\"{o}pf's lemma, we ascertain that \( e(x) \geq k \cdot d(x) \) for some positive constant \( k \), where \( d(x) = \text{dist}(x, \partial\Omega) \).

We define \( C_e(\bar{\Omega}) \) as the set of functions \( u \in C_0(\bar{\Omega}) \) that satisfy \( -t e \leq u \leq t e \) for some \( t > 0 \). The space \( C_e(\bar{\Omega}) \), endowed with the norm \( \|u\|_e = \inf\{t > 0 : -te \leq u \leq te\} \), is an ordered Banach space with a positive cone \( P_e = \{u \in C_e(\bar{\Omega}) : u(x) > 0\} \), which is normal and has a nonempty interior. Specifically, the interior \( \overset{o}{P}_e \) includes all functions \( u \in C_0(\bar{\Omega}) \) bounded by \( t_1 e \leq u \leq t_2 e \) for some \( t_1, t_2 > 0 \). The space relationships are outlined as follows:
\begin{equation}\label{inclusion}
C_0^1(\bar{\Omega}) \hookrightarrow C_e(\bar{\Omega}) \hookrightarrow C_0(\bar{\Omega}).
\end{equation}
Based on the above embedding \eqref{inclusion} and Remark \ref{rem11}, the operator \( S_f: C_e(\bar{\Omega}) \to C_e(\bar{\Omega}) \) is characterized as compact and increasing. Next we will show $S_f$ is strongly increasing.

\begin{lemma}\label{stri}
The operator $S_f:C_{e}(\bar{\Omega})\to C_{e}(\bar{\Omega})$ is strongly increasing.
\end{lemma}
\begin{proof}
Let $u_1,u_2\in C_{e}(\overline{\Omega})$ such that $u_1\leq u_2$ and $u_1\not=u_2$. Define $w_i=S_f(u_i)$ for $i=1,2$, meaning:
\begin{equation}\label{conv22}
\left\{
\begin{aligned}
-\mathcal{M}_{\lambda,\Lambda}^+(D^2w_i)-\Gamma|Dw_i|^2+\gamma w_i&=f(u_i)+\gamma u_i &&\text{in } \Omega,\\
w_i&=0 &&\text{on } \partial\Omega,
\end{aligned}
\right.
\end{equation}
for $i=1,2.$ Since $f(u_1)+\gamma u_1\leq f(u_2)+\gamma u_2$, the comparison principle implies $w_1\leq w_2$. Moreover, the $C^{1,\alpha}(\overline{\Omega})$ estimate (Theorem \ref{c1alpha}) ensures $w_i\in C^{1,\alpha}(\overline{\Omega})$. Rewriting Equation \ref{conv22}, we find that $w_i$ satisfies:
\begin{equation}\label{conv23}
\left\{
\begin{aligned}
-\mathcal{M}_{\lambda,\Lambda}^+(D^2w_i)+\gamma w_i&=F_i(x) &&\text{in } \Omega,\\
w_i&=0 &&\text{on } \partial\Omega,
\end{aligned}
\right.
\end{equation}
where $F_i(x)=f(u_i)+\gamma u_i+\Gamma|Dw_i|^2\in C(\overline{\Omega})$. Therefore, by applying the $W^{2,p}$ estimate, we find that $w_i\in W^{2,p}(\Omega)$ for all $p<\infty$, implying that $w_i$ for $i=1,2$ are strong solutions. From this point, we can treat $w_i$ as if they were classical solutions, and the resulting equations and inequalities will hold in the $L^p$-viscosity sense. Define $w=w_2-w_1$, then $w\geq0$ and satisfies:
\begin{equation}\label{conv23}
\left\{
\begin{aligned}
\mathcal{M}_{\lambda,\Lambda}^-(D^2w)-\gamma b(x)|Dw|-\gamma w&\leq 0 &&\text{in } \Omega,\\
w&\geq 0 &&\text{on } \overline{\Omega},
\end{aligned}
\right.
\end{equation}
where $b(x)=(|Dw_1|+|Dw_2|)\in C(\overline{\Omega})$ and $b\geq0$. By the strong maximum principle (SMP), $w>0$ in $\Omega$. If there exists $x_0\in\partial\Omega$ with $w(x_0) = 0$, the Hopf lemma implies that $\partial_\nu w(x_0)>0$, then $w>0$ in $C_{e}(\bar{\Omega})$. It is important to note that $w_{1},w_{2}\in C^{1,\alpha}(\bar{\Omega})$ and so is $w$, and $w>0$ on $\Omega$.

If it were not possible to find a $t>0$ such that $w>te$, then for each $n\in \mathbf{N}$, there would exist $x_n\in \bar{\Omega}$ such that $w(x_n)<\frac{1}{n}e(x_n)$, consequently, $w(x_n)\leq \frac{1}{n}\|e\|_{L^{\infty}(\Omega)}$. Since $w>0$ in $\Omega$, $x_{n}\rightarrow x$ for some $x\in\partial\Omega$. Consider:
\begin{align*}
\displaystyle{\lim_{n\to \infty}}\left|\frac{w(x_n)-w(x)}{x_n-x}\right|&\leq\displaystyle{\lim_{n\to \infty}}\left|\frac{\frac{1}{n}e(x_n)-\frac{1}{n}e(x)}{x_n-x}\right| \quad (\text{since} \quad w(x)=e(x)=0)\\
&\leq\displaystyle{\lim_{n\to \infty}}\frac{1}{n}\frac{|\nabla e(z_n)||(x_n-x)|}{|x_n-x|}\\
&\leq\displaystyle{\lim_{n\to \infty}}\frac{1}{n}\|\nabla e\|_{L^{\infty}(\Omega)}\to 0,
\end{align*}
This contradicts the fact that \(\frac{\partial w}{\partial \eta} < 0\) on \(\partial \Omega\). Therefore, it must be established that there exists some \( t > 0 \) such that \( w > t e \), confirming that \( S_f \) is indeed a strongly increasing operator.
\end{proof}
\begin{theorem}\label{thm3.3}
The operator \( S_f: C_e(\bar{\Omega}) \to C_e(\bar{\Omega}) \) is completely continuous.
\end{theorem}

\begin{proof}
It is already established that \( S_f \) is compact. We now demonstrate its continuity. Consider a sequence \( u_n \in C_e(\bar{\Omega}) \) that converges to \( u \) in \( C_e(\bar{\Omega}) \). According to Remark \ref{rem11}, this convergence in \( C_e(\bar{\Omega}) \) ensures that \( S_f(u_n) \) converges to \( S_f(u) \) in \( C^1(\bar{\Omega}) \). Furthermore, given the inclusion property \eqref{inclusion}, it follows that \( S_f(u_n) \) converges to \( S_f(u) \) in \( C_e(\Omega) \), thus confirming the complete continuity of \( S_f \).
\end{proof}
\section{Main Results}
\begin{lemma}\label{fixthm}
A function \( u \) solves \eqref{problem} if and only if \( u \) is a fixed point of the map \( S_f \) in \( C_{e} \).
\end{lemma}
\begin{proof}
Suppose \( u \) is a solution to \eqref{problem}. According to Theorem \ref{c1alpha} and the embedding \eqref{inclusion}, \( u \) belongs to \( C_{e} \) and \( S_f(u) = u \) by the definition of \( S_f \). Conversely, if \( u \) in \( C_{e} \) is a fixed point of \( S_f \), Theorem \ref{c1alpha} ensures \( u \) belongs to \( C^{1,\alpha}_{0}(\bar{\Omega}) \). Hence, the function \( g(u) = \gamma u + f(u) + \Gamma|Du|^2 \) is H\"{o}lder continuous, hence \( u \in C^2 \)  and  \( u \) is a solution of \eqref{problem}. \end{proof}

Thus, the establishment of three solutions to \eqref{problem} is equivalent to demonstrating the existence of three fixed points for \( S_f \). We plan to validate this by applying Lemma \ref{three}, as stated below. For this application, it is essential to confirm that \( S_f \) fulfills all necessary conditions.
\begin{theorem}{(Minimal and Maximal Solutions)}\label{Th3.5}
Assume $\psi$ and $\phi$ are positive sub and supersolutions of \eqref{problem}, respectively, satisfying $\psi \leq \phi$. Then, there exist minimal and maximal solutions to \eqref{problem} within the ordered interval $[\psi, \phi]$.
\end{theorem}
\begin{proof}
Given that $\psi$ is a subsolution of \eqref{problem}, then:
\begin{align*}
-\mathcal{M}_{\lambda,\Lambda}^+(D^2\psi) - \Gamma|D\psi|^2 + \gamma \psi &\leq f(\psi) + \gamma \psi\\
&= -\mathcal{M}_{\lambda,\Lambda}^+(D^2 S_f(\psi)) - \Gamma|D S_f(\psi)|^2 + \gamma S_f(\psi).
\end{align*}
By the comparison principle, it follows that $\psi \leq S_f(\psi)$. Similarly, $S_f(\phi) \leq \phi$ and $S_f(\psi) \leq S_f(\phi)$ (due to the monotonicity of $S_f$). Let $E = C_e(\bar{\Omega})$ and let $P$ be the positive cone of $E$. Given that $E$ is an ordered Banach space, Theorem \ref{stri} and \ref{thm3.3} imply that $S_f : [\psi, \phi] \rightarrow E$ is both increasing and completely continuous. Consequently, $\psi \leq S_f(\psi) \leq S_f(\phi) \leq \phi$. Therefore, by Corollary 6.2\cite{aman}, there exist minimal and maximal fixed points of $S_f$ within the interval $[\psi, \phi]$. According to Lemma \ref{fixthm}, these fixed points correspond to the minimal and maximal solutions of \eqref{problem} within $[\psi, \phi]$.
\end{proof}

\begin{lemma}\label{three}[\cite{aman}]
Suppose $X$ is a retract of a Banach space and $F : X \to X$ is a completely continuous map. If $X_1$ and $X_2$ are disjoint retracts within $X,$ and let $U_k, k=1,2$ be open subsets within $X_k, k=1,2$. Additionally, assume $F(X_k) \subset X_k$ and $F$ has no fixed points in $X_k \setminus U_k, k=1,2$. Then $F$ will have at least three distinct fixed points, namely $x, x_1, x_2$ with $x_k \in X_k, k=1,2$ and $x \in X \setminus (X_1 \cup X_2)$.
\end{lemma}

\begin{theorem}[Three Solution Theorem]
\label{Th3.7}
Assume the existence of two pairs of ordered sub and supersolutions, $(\psi_1, \phi_1)$ and $(\psi_2, \phi_2)$, for \eqref{problem} satisfying:
\begin{itemize}
    \item $\psi_1 \leq \psi_2 \leq \phi_1$ and $\psi_1 \leq \phi_2 \leq \phi_1$,
    \item $\psi_2 \not\leq \phi_2$, and
    \item neither $\psi_2$ nor $\phi_2$ are solutions of \eqref{problem}.
\end{itemize}
Then, there are at least three solutions $u_1, u_2, u_3$ to \eqref{problem} such that:
\begin{itemize}
    \item $u_1 \in [\psi_1, \phi_2]$,
    \item $u_2 \in [\psi_2, \phi_1]$, and
    \item $u_3 \in [\psi_1, \phi_1] \setminus ([\psi_1, \phi_2] \cup [\psi_2, \phi_1])$.
\end{itemize}
\end{theorem}

\begin{proof}
Invoking Lemma \ref{fixthm}, we identify solutions to \eqref{problem} with fixed points of the operator $S_f$. Let us define the sets $X = [\psi_1, \phi_1]$, $X_1 = [\psi_1, \phi_2]$, and $X_2 = [\psi_2, \phi_1]$. The inequalities among $\phi's$ and $\psi's$ ensure that $X_1, X_2 \subset X$ and $X_1 \cap X_2 = \emptyset$. Since $S_f$ is an increasing operator and both $\psi_1, \phi_1$ are ordered sub and supersolutions, it follows that:
\begin{align*}
    S_f(X) \subseteq X, \quad S_f(X_1) \subseteq X_1, \quad S_f(X_2) \subseteq X_2.
\end{align*}
$S_f$ is also a retract of $X$ and $X_k$ for $k = 1, 2$ in the Banach space $C_e(\bar{\Omega})$. It is to see from Lemma \ref{stri} and Theorem \ref{thm3.3} that $S_f$ is strongly increasing and completely continuous when restricted to $X$, $X_1$, and $X_2$.

It is given that $\phi_2$ is a strict supersolution, meaning:
\begin{equation}
-\mathcal{M}_{\lambda,\Lambda}^+(D^2 \phi_2) - \Gamma |D\phi_2|^2 + \gamma \phi_2 \geq f(\phi_2) + \gamma \phi_2 + g,
\end{equation}
for some positive function $g$ on $\overline{\Omega}$. The operator definition of $S_f$ leads to:
\begin{equation}
-\mathcal{M}_{\lambda,\Lambda}^+(D^2S_f(\phi_2)) - \Gamma|DS_f(\phi_2)|^2 + \gamma S_f(\phi_2) = f(\phi_2) + \gamma \phi_2.
\end{equation}
By the viscosity solution concept, this implies:
\begin{equation}
-\mathcal{M}_{\lambda,\Lambda}^+(D^2(S_f(\phi_2) - \phi_2)) - \Gamma \tilde{b}|D(S_f(\phi_2) - \phi_2)| + \gamma (S_f(\phi_2) - \phi_2) \leq -g < 0,
\end{equation}
where $\tilde{b} = |DS_f(\phi_2)| + |D\phi_2|$. The comparison principle confirms $S_f(\phi_2) < \phi_2$ in $\Omega$. By the strong maximum principle, we conclude $S_f(\phi_2) \neq \phi_2$. Thus, $S_f$ has a maximal fixed point $u_1 \in X_1$ with $\psi_1 \leq u_1 < \phi_2$ and a minimal fixed point $u_2 \in X_2$ with $\psi_2 < u_2 \leq \phi_1$.

The existence of a third fixed point $u_3$, distinct from $u_1$ and $u_2$, follows by Lemma \ref{three} which ensures that $u_3$ must be located in $[\psi_1, \phi_1] \setminus ([\psi_1, \phi_2] \cup [\psi_2, \phi_1])$.
\end{proof}

\begin{proof}
Given that $\phi_2$ is designated as a strict supersolution of \eqref{problem}, we have:
\begin{equation}\label{ww}
-\mathcal{M}_{\lambda,\Lambda}^+(D^2\phi_2) - \Gamma |D\phi_2|^2 + \gamma \phi_2 \geq f(\phi_2) + \gamma \phi_2 + g,
\end{equation}
where $g \in C(\bar{\Omega})$ is positive throughout $\bar{\Omega}$. Additionally, the operator definition for $S_f$ provides:
\begin{equation}\label{ww2}
-\mathcal{M}_{\lambda,\Lambda}^+(D^2 S_f(\phi_2)) - \Gamma |D S_f(\phi_2)|^2 + \gamma S_f(\phi_2) = f(\phi_2) + \gamma \phi_2.
\end{equation}
Comparing \eqref{ww} and \eqref{ww2} yields the viscosity inequality:
\begin{equation}\label{super}
-\mathcal{M}_{\lambda,\Lambda}^+(D^2(S_f(\phi_2) - \phi_2)) - \Gamma \tilde{b} |D(S_f(\phi_2) - \phi_2)| + \gamma (S_f(\phi_2) - \phi_2) \leq -g < 0,
\end{equation}
where $\tilde{b} = |D S_f(\phi_2)| + |D \phi_2|$. The comparison principle then implies $S_f(\phi_2) \leq \phi_1$ in $\Omega$. By the strong maximum principle, we deduce either $S_f(\phi_2) = \phi_2$ or $S_f(\phi_2) < \phi_2$. Given that $\phi_2$ is a strict supersolution, $S_f(\phi_2) \equiv \phi_2$ is untenable, leading to $S_f(\phi_2) < \phi_2$. Consequently, by Corollary 6.2\cite{aman}, $S_f$ admits a maximal fixed point $u_1$ in $X_1$, satisfying $\psi_1 \leq u_1 < \phi_2$. Similarly, as $\psi_2$ is a strict subsolution, $S_f$ possesses a minimal fixed point $u_2$ in $X_2$, where $\psi_2 < u_2 \leq \phi_1$.

For $u_1$ being a solution and $\phi_2$ a strict supersolution, the inequality holds:
\begin{equation}
-\mathcal{M}_{\lambda,\Lambda}^-(D^2(\phi_2 - u_1)) - \Gamma \tilde{b} |D(\phi_2 - u_1)| + \gamma (\phi_2 - u_1) \geq f(\phi_2) - f(u_1) + \gamma (\phi_2 - u_1) + g,
\end{equation}
where $\tilde{b} = |D \phi_2| + |D u_1|$. As $f(s) + \gamma s$ is an increasing function and $u_1 \leq \phi_2$, the above inequality simplifies to:
\begin{equation}
-\mathcal{M}_{\lambda,\Lambda}^-(D^2(\phi_2 - u_1)) - \Gamma \tilde{b} |D(\phi_2 - u_1)| + \gamma (\phi_2 - u_1) \geq g > 0.
\end{equation}
Reiterating the arguments in Theorem \ref{stri}, constants $t_1 > 0$ and $t_2 > 0$ are found such that $\phi_2 - u > t_1 e$ and $u_2 - \psi_2 > t_2 e$. Defining:
\begin{equation}
B_k = X \cup \{z \in C_e(\bar{\Omega}) : \|z - u_k\| < t_k\} \quad \text{for } k = 1, 2,
\end{equation}
ensures $B_k \subset X_k$ and each $B_k$ is open in $X$, where $X_k$ has a non-empty interior. Let $U_k$ be the largest open set in $X_k$ containing $u_k$ with $S_f$ having no fixed points in $X_k \setminus U_k$ for $k=1, 2$. The existence of a third fixed point for $S_f$ follows from Lemma \ref{Th3.7}, thereby affirming that \eqref{problem} admits at least three solutions: $u_1 \in [\psi_1, \phi_2]$, $u_2 \in [\psi_2, \phi_1]$, and $u_3 \in [\psi_1, \phi_1] \setminus ([\psi_1, \phi_2] \cap [\psi_2, \phi_1])$.
\end{proof}

\section{Results Applied to the Sublinear Case}
Building upon the foundation established in the previous section, we now demonstrate the existence of three distinct solutions to the boundary value problem defined as:
\begin{equation}\label{prob}
\left\{
\begin{aligned}
-\mathcal{M}_{\lambda,\Lambda}(D^2 u) - \Gamma |Du|^2 &= \mu g(u) && \text{in } \Omega, \\
u &= 0 && \text{on } \partial \Omega,
\end{aligned}
\right.
\end{equation}
where \( g: [0, \infty) \to \mathbb{R} \) is a nondecreasing H\"older continuous function with \( g(0) = 0 \). This problem's solution strategy leverages Theorem \ref{Th3.7} to affirm the existence of multiple positive solutions.

Key to applying Theorem \ref{Th3.7} is the identification and construction of an ordered pair of sub and supersolutions, a technique adapted from the methodologies discussed in \cite{ramaswamy2004multiple, Mallick_Verma_2021}. The principal challenge here involves the superlinear growth condition in the gradient term, which significantly complicates the problem.

To address this issue, we employ a transformation of the dependent variable, thoroughly explored and justified in Lemma 2.3 of \cite{sirakov2010}. This transformation effectively mitigates the complexities introduced by the gradient's superlinear growth, thereby allowing the subsequent application of Theorem \ref{Th3.7}.

\begin{lemma}\label{lemma1}
Let \( u \in W_{\text{loc}}^{2,N}(\Omega) \). For any positive constant \( m \), define the transformation \( v = \frac{e^{mu} - 1}{m} \). Then, almost everywhere in \( \Omega \), the following gradient relationship holds:
\[ Dv = (1 + mv)Du, \]
and the inequalities are given by
\[ m\lambda|Du|^2 + \mathcal{M}_{\lambda,\Lambda}^{+}(D^2u) \leq \frac{\mathcal{M}_{\lambda,\Lambda}^{+}(D^2v)}{1 + mv} \leq m\Lambda|Du|^2 + \mathcal{M}_{\lambda,\Lambda}^{+}(D^2u). \]
Furthermore, if \( u = 0 \) on \( \partial\Omega \) and \( u > 0 \) in \( \Omega \), then \( v \) also satisfies \( v = 0 \) on \( \partial\Omega \) and \( v > 0 \) within \( \Omega \).
\end{lemma}
This lemma holds true for both strong and viscosity solutions. The transformation introduced here effectively handles the challenges posed by the quadratic gradient term but introduces complexity in managing the nonlinear terms.

To establish the existence of a positive solution to \eqref{prob} when \( \mu \) is small, we will demonstrate that as \( \mu \) approaches zero, the norm \( \|u_{\mu}\|_{\infty} \) also tends towards zero. This proof relies on specific assumptions about the nondecreasing function \( g \), which are as follows:

\begin{itemize}
\item[\textbf{(C1)}] $\displaystyle{\lim_{s \to 0} \frac{g(s)}{s} = \infty}.$
\end{itemize}

\begin{theorem}\label{Th4.1}
Assuming conditions \(C0\) and \(C1\), there exists a \(\mu_0 > 0\) such that for all \(\mu\) in the range \(0 \leq \mu \leq \mu_0\), the problem \eqref{prob} admits a positive solution \(u \in C^2(\Omega) \cap C(\bar{\Omega})\). Furthermore, the norm \(\|u_{\mu}\|_{\infty}\) approaches zero as \(\mu\) tends to $0$.
\end{theorem}

\begin{proof}
Start by selecting \(e \in C^2(\Omega) \cap C(\bar{\Omega})\), as described in \eqref{ee}. Given that \(g(0) = 0\) and \(g\) is an increasing function, we select a sufficiently small \(\mu_0 > 0\) such that
\[
(1 + m_1\mu_0\|e\|_{\infty}) g\left(\frac{1}{m_1} \log(1 + m_1\mu_0\|e\|_{\infty})\right) < 1,
\]
where \(m_1 = \frac{\Gamma}{\Lambda}\). For \(\mu \leq \mu_0\), we define \(\phi = \mu e\). The nondecreasing property of \(g\) leads to:
\begin{align*}
-\mathcal{M}_{\lambda,\Lambda}^+(D^2\phi) &\geq \mu (1 + m_1\mu e) g\left(\frac{1}{m_1} \log(1 + m_1\mu e)\right) \\
&= \mu (1 + m_1\phi) g\left(\frac{1}{m_1} \log(1 + m_1\phi)\right).
\end{align*}
From Lemma \ref{lemma1}, it follows that \(\tilde{\phi} = \frac{1}{m_1} \log(1 + m_1\phi)\) acts as a supersolution to \eqref{prob} for all \(\mu \leq \mu_0\).

Next, to construct a subsolution, let \(\mu^{+}_{1}\) and a function \(\phi^{+}_{1} \in C^2(\Omega) \cap C(\bar{\Omega})\), satisfying:
\begin{equation}\label{eigenvalue}
-\mathcal{M}^{+}_{\lambda,\Lambda}(D^2\phi^{+}_{1}) = \mu^{+}_{1} \phi^{+}_{1} \quad \text{in } \Omega, \quad \phi^{+}_{1} = 0 \text{ on } \partial\Omega,
\end{equation}
with \(\phi^{+}_{1} > 0\) in \(\Omega\) and \(\|\phi_1^+\|_{\infty} = 1\). where \(\mu^{+}_{1}\) and \(\phi^{+}_{1}\) are the first half eigenvalue and eigenfunction respectively to the above problem (see \cite{MR2124162} Proposition 1.1).

Using \(C1\), for a fixed \(\mu > 0\), we find \(m_{\mu} > 0\) such that:
\[
\frac{\mu_1^+}{\mu} \leq \frac{g\left(\frac{1}{m_2} \log(1 + m_2 m_{\mu} \phi_1^+)\right)}{m_{\mu} \phi_1^+},
\]
where \(m_2 = \frac{\Gamma}{\lambda}\). Letting \(\psi = m_{\mu} \phi_1^+\), it follows that:
\[
\mathcal{M}_{\lambda,\Lambda}^+(D^2\psi) \leq \mu g\left(\frac{1}{m_2} \log(1 + m_2 m_{\mu} \phi_1^+)\right)\leq (1 + m_2 m_{\mu} \phi_1^+)g\left(\frac{1}{m_2} \log(1 + m_2 m_{\mu} \phi_1^+)\right).
\]
Thus, \(\tilde{\psi} = \frac{1}{m_2} \log(1 + m_2 \psi)\) becomes a subsolution for \eqref{prob}. Since $\frac{\partial e}{\partial \eta} < 0$ on $\partial\Omega$, hence $\frac{\partial \tilde{\phi}}{\partial \eta} = \frac{\mu}{1 + m_1 \phi} \frac{\partial e}{\partial \eta} < 0$, by choosing $m_{\mu}$ sufficiently small, we have \(\tilde{\psi} \leq \tilde{\phi}\), Theorem \ref{Th3.5} guarantees the existence of a solution \(u_{\mu}\) within the bounds \(\tilde{\psi} \leq u_{\mu} \leq \tilde{\phi}\) for \(0 \leq \mu \leq \mu_0\). Moreover, \(\|u_{\mu}\|_{\infty}\) approaches zero as \(\mu \to 0\), affirming that the norm of the solution shrinks with diminishing \(\mu\).
\end{proof}

To further this discussion, we assume an additional condition on \(g\):
\begin{itemize}
\item[\textbf{(C2)}] \(\displaystyle{\lim_{s \to \infty} g(s) = c}\),
\end{itemize}
where \(c > 0\). In the next theorem that will establish the existence of a positive solution to \eqref{prob} under the new constraints for certain \(\mu > 0\).
\begin{theorem}\label{Th4.2}
Assuming conditions (C0) through (C2), there exists \(\hat{\mu} > 0\) dependent on the ellipticity constant \(m_1 = \frac{\Gamma}{\Lambda}\) and the limit \(c\). For all \(\mu\) in the range \(0 \leq \mu \leq \hat{\mu}\), problem \eqref{prob} has a positive solution \(u_{\mu}\).
\end{theorem}
\begin{proof}
Recall from the previous theorem, \(\tilde{\psi} = \frac{1}{m_2} \log(1 + m_2 \psi)\) with \(\psi = m_{\mu} \phi_1^+\) acts as a valid subsolution. We now construct a supersolution \(\tilde{\phi}\) for problem \eqref{prob}. Define \(\hat{\mu} = \frac{1}{c(m_1 + 1) \|e\|_{\infty}}\), and for any \(\mu\) less than \(\hat{\mu}\), choose \(M_{\mu} > 0\) such that \(M_{\mu} \|e\|_{\infty} > 1\). Given that \(\mu \leq \hat{\mu}\), it follows that:
\[
\frac{1}{\mu \|e\|_{\infty}} \geq c(1 + m_1),
\]
implying:
\[
\left(\frac{1}{M_{\mu} \|e\|_{\infty}} + m_1\right) g\left(\frac{1}{m_1} \log(1 + m_1 M_{\mu} \|e\|_{\infty})\right) \leq c(1 + m_1).
\]
From this, we deduce:
\[
M_{\mu} \geq \mu (1 + m_1 M_{\mu} \|e\|_{\infty}) g\left(\frac{1}{m_1} \log(1 + m_1 M_{\mu} \|e\|_{\infty})\right).
\]
Defining \(\phi = M_{\mu} e\), we obtain:
\[
-\mathcal{M}_{\lambda,\Lambda}^+(D^2\phi) \geq \mu (1 + m_1 \phi) g\left(\frac{1}{m_1} \log(1 + m_1 \phi)\right),
\]
confirming that \(\tilde{\phi} = \frac{1}{m_1} \log(1 + m_1 \phi)\) serves as a supersolution to \eqref{prob} for all \(\mu \leq \hat{\mu}\). Given \(\frac{\partial e}{\partial \eta} < 0\) on \(\partial\Omega\) and \(M_{\mu} \|e\|_{\infty} > 1\), we ensure \(m_{\mu}\) is sufficiently small such that \(\tilde{\psi} \leq \tilde{\phi}\). Hence, for each \(\mu \in (0 , \hat{\mu})\), \eqref{prob} admits a solution \(u_{\mu}\) satisfying \(\tilde{\psi} \leq u_{\mu} \leq \tilde{\phi}\).
\end{proof}
The next step will investigate conditions under which \eqref{prob} fails to have positive solutions when \(\mu\) is large.

\begin{theorem}\label{nonexist}
There exists $\tilde{\mu} > \mu_1^+$ dependent on $c$ and $m$ such that \eqref{prob} does not have a positive solution for $\mu > \tilde{\mu}$.
\end{theorem}
\begin{proof}
Since $g$ is non-decreasing, $g'(0) > 0$, and $g$ is bounded above by $c$, let us write $f(s) = (1 + m_1 s) g\left(\frac{1}{m_1} \log(1 + m_1 s)\right)$, where $m_1 = \frac{\Gamma}{\Lambda}$. Hence, $\displaystyle{\lim_{s \to \infty} \frac{f(s)}{s} = m_1 c > 0}$, so there exists $s_0 > 0$ such that:
$$\left|\frac{f(s)}{s} - m_1 c\right| < \frac{m_1 c}{2} \text{ for all } s > s_0.$$
Hence, $\frac{m_1 c}{2} \leq \frac{f(s)}{s} \leq \frac{3 m_1 c}{2}$ for $s > s_0$. As $f'(0) > 0$ and $f$ is non-decreasing, it is easy to see that there exists $k > 0$ such that $f(s) > ks$ for all $s \in [0, s_0]$, where $k = \inf\{l > 0 : f(s) > ls, \text{ for all } s \in [0, s_0]\}$. Hence, $f(s) > \alpha s$ for all $s > 0$, where $\alpha = \min\{k, \frac{3 m_1 c}{2}\}$.

Let us assume that \eqref{prob} has a positive solution $u$ for $\mu > \frac{\mu_1^+}{\alpha} = \tilde{\mu}$, i.e.:
\begin{equation*}
\left\{
\begin{aligned}
-\mathcal{M}_{\lambda,\Lambda}(D^2u)-\Gamma|Du|^2 &= \mu f(u) &&\text{in } \Omega,\\
u &= 0 &&\text{on } \partial\Omega,
\end{aligned}
\right.
\end{equation*}
Define $v = \frac{e^{m_1 u} - 1}{m_1}$. Then from Lemma \ref{lemma1}, we have:
$$-\mathcal{M}^+_{\lambda,\Lambda}(D^2 v) \geq \mu (1 + m_1 \phi) g\left(\frac{1}{m_1} \log(1 + m_1 \phi)\right) \geq \mu \alpha v.$$
Let $\phi_1^+$ be the first half eigenfunction with eigenvalue $\mu_1^+$ of:
\begin{equation*}
\left\{
\begin{aligned}
-\mathcal{M}^+_{\lambda,\Lambda}(D^2\phi_1^+) &= \mu_1^+ \phi_1^+ &&\text{in } \Omega,\\
\phi_1^+ &= 0 &&\text{on } \partial\Omega,
\end{aligned}
\right.
\end{equation*}
then $-\mathcal{M}^+_{\lambda,\Lambda}(D^2\phi_1^+) = \mu_1^+ \phi_1^+ < \mu \alpha \phi_1^+$. Hence, we conclude that for $\mu > \tilde{\mu}$, the problem:
\begin{equation*}
\left\{
\begin{aligned}
-\mathcal{M}^+_{\lambda,\Lambda}(D^2 u) &= \mu \alpha u &&\text{in } \Omega,\\
u &= 0 &&\text{on } \partial\Omega,
\end{aligned}
\right.
\end{equation*}
has a positive solution. This contradicts the maximality of $\mu^{+}_1$ in Proposition 1.1(iv) of \cite{esteban2005nonlinear}.
\end{proof}
The forthcoming theorem explores the conditions under which multiple solutions to Problem \eqref{prob} exist within a specific range of the parameter $\mu$. To establish this theorem, it is necessary to introduce an additional condition on the function $g$:
\begin{itemize}
\item[\textbf{(C3)}] Suppose there exist two positive constants $0 < a < b$ with $g(a) \neq 0$ and $g(b) \neq 0$ such that $\displaystyle{\min\left(\hat{\mu}, \frac{a}{\|e\|_{\infty}(1 + m_1 a)g(a)}\right) > A \frac{b}{g\left(\frac{1}{m_2} \log(1 + m_2 b)\right)},}$
\end{itemize}
where $A = \displaystyle{\inf_{\epsilon} \frac{\Lambda N_{-} R^{N_+ - 1}}{\epsilon^{N_{-}} (R - \epsilon)}}$ and $R$ is the radius of the largest inscribed ball $B_R$ in $\Omega$.

\begin{theorem}\label{nonexist}
Assuming conditions \(C0\) through \(C2\), there exists a critical value \(\tilde{\mu} > \mu_1^+\) that depends on the upper bound \(c\) and the ellipticity ratio \(m_1 = \frac{\Gamma}{\Lambda}\). For any \(\mu > \tilde{\mu}\), problem \eqref{prob} does not admit a positive solution.
\end{theorem}

\begin{proof}
Given that \(g\) is non-decreasing with \(g'(0) > 0\) and bounded above by \(c\), define the function 
\[ f(s) = (1 + m_1 s) g\left(\frac{1}{m_1} \log(1 + m_1 s)\right). \]
This function satisfies
\[ \lim_{s \to \infty} \frac{f(s)}{s} = m_1 c > 0, \]
implying there exists \(s_0 > 0\) such that 
\[ \left|\frac{f(s)}{s} - m_1 c\right| < \frac{m_1 c}{2} \quad \text{for all } s > s_0. \]
Thus, for \(s > s_0\), 
\[ \frac{m_1 c}{2} \leq \frac{f(s)}{s} \leq \frac{3 m_1 c}{2}. \]
For \(s\) in the range \([0, s_0]\), \(f(s)\) exceeds \(ks\), where \(k\) is the lowest bound ensuring 
\[ f(s) > ks, \quad k = \inf\{l > 0 : f(s) > ls \text{ for all } s \in [0, s_0]\}. \]
This is possible as $f$ is nondeareasing and $f'(0)>0$. Thus, for all \(s > 0\),
\[ f(s) > \alpha s, \quad \alpha = \min\{k, \frac{3 m_1 c}{2}\}. \]
Assuming the existence of a positive solution \(u\) to \eqref{prob} for \(\mu > \frac{\mu_1^+}{\alpha} = \tilde{\mu}\):
\begin{align*}
-\mathcal{M}^+_{\lambda,\Lambda}(D^2u)-\Gamma|Du|^2 &= \mu f(u) &&\text{in } \Omega,\\
u &= 0 &&\text{on } \partial\Omega,
\end{align*}
define \(v = \frac{e^{m_1 u} - 1}{m_1}\). From Lemma \ref{lemma1}, we know
\[ -\mathcal{M}^+_{\lambda,\Lambda}(D^2 v) \geq \mu \alpha v. \]
Let \(\phi_1^+\) be the principal half eigenfunction associated with eigenvalue \(\mu_1^+\) of \ref{eigenvalue}. Since
\[ -\mathcal{M}^+_{\lambda,\Lambda}(D^2\phi_1^+) = \mu_1^+ \phi_1^+ < \mu \alpha \phi_1^+, \]
for \(\mu > \tilde{\mu}\), this contradicts the maximality of \(\mu_1^+\) established in Proposition 1.1(iv) from \cite{esteban2005nonlinear}, implying that no positive solutions exist for \(\mu > \tilde{\mu}\).
\end{proof}

In the next theorem, we will explore the existence of multiple solutions to Problem \eqref{prob} over a specific range of \(\mu\) values. For this proof, we impose an additional condition on \(g\):

\begin{itemize}
\item[\textbf{(C3)}] Assume two positive constants \(0 < a < b\) with \(g(a) \neq 0\) and \(g(b) \neq 0\) are such that
\[ \min\left(\hat{\mu}, \frac{a}{\|e\|_{\infty}(1 + m_1 a)g(a)}\right) > A \frac{b}{g\left(\frac{1}{m_2} \log(1 + m_2 b)\right)}, \]
where \(A = \displaystyle{\inf_{\epsilon} \frac{\Lambda N_{-} R^{N_+ - 1}}{\epsilon^{N_{-}} (R - \epsilon)}}\) and \(R\) is the radius of the largest inscribed ball \(B_R\) in \(\Omega\).
\end{itemize}

\begin{theorem}\label{Th4.3}
Let \( g \) satisfy conditions \( C0 \) through \( C3 \). Then there exist \( 0 \leq \mu_* < \mu^* \) such that for any \( \mu \) satisfying \( \mu_* < \mu < \mu^*\), Problem \eqref{prob} has three positive solutions.
\end{theorem}

\begin{proof}
To prove the existence of three solutions, we apply Theorem \ref{Th3.7}, which necessitates constructing appropriate pairs of sub and supersolutions. These pairs are denoted \( \tilde{\psi}_1, \tilde{\psi}_2 \) and \( \tilde{\phi}_1, \tilde{\phi}_2 \), respectively, which need to satisfy the conditions outlined in Theorem \ref{Th3.7}. The pairs \( \tilde{\psi}_1 \) and \( \tilde{\phi}_1 \) can be derived from Theorem \ref{Th4.2}.

For \( \tilde{\psi}_2 \) and \( \tilde{\phi}_2 \), consider \( \Omega \) as the ball \( B_R \) initially, to simplify the construction process. Define \( \phi_2 = \frac{a e}{\|e\|_\infty} \), where $e \in C^2(\Omega) \cap C(\bar{\Omega})$ from Theorem \ref{Pr2.6}. Thus, for $\mu < \frac{a}{\|e\|_\infty (1 + m_1 a)g(a)}$, we have:
\begin{align*}
-\mathcal{M}_{\lambda,\Lambda}^+ (D^2 \phi_2) &= \frac{a}{\|e\|_\infty} > \mu (1 + m_1 a) g(a) \\
&> \mu (1 + m_1 a) g\left(\frac{1}{m_1} \log(1 + m_1 a)\right) \\
&> \mu \left(1 + m_1 \frac{a e}{\|e\|_{\infty}}\right) g\left(\frac{1}{m_1} \log\left(1 + m_1 \frac{a e}{\|e\|_{\infty}}\right)\right) \\
&= \mu (1 + m_1 \phi_2) g\left(\frac{1}{m_1} \log(1 + m_1 \phi_2)\right).
\end{align*}
Therefore, \( \tilde{\phi}_2 = \frac{1}{m_1} \log(1 + m_1 \phi_2) \) serves as a strict supersolution for \( \mu \leq \frac{a}{\|e\|_\infty (1 + m_1 a) g(a)} \).

Next our aim to construct a positive strict subsolution $\tilde{\psi}_2$ for problem \eqref{prob}. Define $\mu > \mu_* = A \frac{b}{g\left(\frac{1}{m_2} \log(1 + m_2 b)\right)}$ and let $0 < \epsilon < R$ and $l, m > 1$. Introduce a function $\rho: [0, R] \to [0, 1]$ by
\[
\rho(r) = 
\begin{cases}
1 & \text{if } 0 \leq r \leq \epsilon,\\
1 - \left(1 - \left(\frac{R - r}{R - \epsilon}\right)^m\right)^l & \text{if } \epsilon < r \leq R.
\end{cases}
\]
Then the derivative of $\rho$, $\rho'(r)$, is given by
\[
\rho'(r) = 
\begin{cases}
0 & \text{if } 0 \leq r \leq \epsilon,\\
-\frac{lm}{R - \epsilon} \left(1 - \left(\frac{R - r}{R - \epsilon}\right)^m\right)^{l - 1} \left(\frac{R - r}{R - \epsilon}\right)^{m - 1} & \text{if } \epsilon < r \leq R.
\end{cases}
\]
Consequently, $|\rho'(r)| \leq \frac{lm}{R - \epsilon}$. Letting $d(r) = b \rho(r)$, we find that
\begin{equation}
|d'(r)| \leq \frac{b lm}{R - \epsilon}. \label{10}
\end{equation}

Consider $\psi_2$ as a solution to
\begin{equation}\label{211}
\begin{aligned}
-\mathcal{M}_{\lambda,\Lambda}^+ (D^2\psi_2) &= \mu g\left(\frac{1}{m_2} \log(1 + m_2 d(r))\right) & \text{in } B_R,\\
\psi_2 &= 0 & \text{on } \partial B_R.
\end{aligned}
\end{equation}
By established regularity results, $\psi_2 \in C^{2}(\overline{B_R})$ and is radially symmetric (Theorem 1.1 in \cite{sym}). Thus, the differential equations simplify to
\begin{equation}\label{12}
\begin{aligned}
-&\theta(\psi_2''(r)) \psi_2''(r) - \theta(\psi_2'(r)) \frac{(N - 1)}{r} \psi_2'(r) = \mu g\left(\frac{1}{m_2} \log(1 + m_2 d(r))\right) & \text{in } B(0, R),\\
&\psi_2'(0) = \psi_2(R) = 0,
\end{aligned}
\end{equation}
where $\theta(s) = \Lambda \text{ if } s > 0 \text{ and } \lambda \text{ if } s \leq 0$.

By transforming \eqref{12} using the variable substitution $v(r) = \exp\left({\int_1^r \frac{\theta(\psi_2'(s))(N - 1)}{\theta(\psi_2''(s)) s} ds}\right)$ and $\tilde{v}(r) = \frac{v(r)}{\theta(\psi_2''(r))}$, we arrive at
\begin{equation*}\label{3.2}
\begin{aligned}
v(r) \psi_2''(r) + v(r) \frac{\theta(\psi_2'(r))}{\theta(\psi_2''(r))} \frac{(N - 1)}{r} \psi_2'(r) &= -\mu \tilde{v}(r) g\left(\frac{1}{m_2} \log(1 + m_2 d(r))\right),
\end{aligned}
\end{equation*}
which implies
\begin{equation}\label{3.21}
\begin{aligned}
\left(v(r) \psi_2'(r)\right)' &= -\mu \tilde{v}(r) g\left(\frac{1}{m_2} \log(1 + m_2 d(r))\right).
\end{aligned}
\end{equation}
From \eqref{3.21}, integrating from $0$ to $r$ shows
\[
v(r) \psi_2'(r) = -\frac{\mu}{v(r)} \int_0^r \tilde{v}(s) g\left(\frac{1}{m_2} \log(1 + m_2 d(s))\right) ds \quad \text{(since } \psi_2'(0) = 0).
\]
We claim that for $\mu > \mu_*$, $\psi_2(t) > d(t)$ for $t \in [0, R)$. Assuming the claim, $-\mathcal{M}_{\lambda,\Lambda}^+(D^2 \psi_2(r))$ will be less than $\mu (1 + m_2 \psi_2(r)) g\left(\frac{1}{m_2} \log(1 + m_2 \psi_2(r))\right)$ in $B_R$, showing that $\psi_2$ is a strict supersolution of \eqref{prob}.
To demonstrate $d(t) < \psi_2(t)$, it suffices to show $\psi_2'(t) < d'(t)$ on $(0, R]$ since $\psi_2(R) = d(R) = 0$. For any $r \in (0, \epsilon)$, $g\left(\frac{1}{m_2} \log(1 + m_2 d(r))\right) = g\left(\frac{1}{m_2} \log(1 + m_2 b)\right)$, and $\tilde{\tau}(r) > 0$, hence $\psi_2'(r) < 0 = d'(r)$, supporting the claim. For $r \in (\epsilon, R)$, let us define:
\begin{align*}
N_+ &= \frac{\lambda}{\Lambda}(N - 1) + 1,\\
N_- &= \frac{\Lambda}{\lambda}(N - 1) + 1,
\end{align*}
and note that:
\begin{itemize}
    \item $N_+ - 1 \leq v(r) r \leq N_- - 1$,
    \item $r^{N_- - 1} \leq \tau(r) \leq r^{N_+ - 1}$ and $\frac{\tau(r)}{\Lambda} \leq \tilde{\tau}(r) \leq \frac{\tau(r)}{\lambda}$.
\end{itemize}

Consider the inequality:
\begin{align*}
-\psi_2'(r) &= \frac{\mu}{\tau(r)} \int_0^r \tilde{\tau}(s) g\left(\frac{1}{m_2} \log(1 + m_2 d(s))\right) ds \\
&\geq \frac{\mu}{r^{N_+ - 1}} \int_0^\epsilon \tilde{\tau}(s) g\left(\frac{1}{m_2} \log(1 + m_2 d(s))\right) ds \\
&\geq \frac{\mu g\left(\frac{1}{m_2} \log(1 + m_2 b)\right) \epsilon^{N_-}}{\Lambda R^{N_+ - 1} N_-}.
\end{align*}

To establish $-\psi_2'(r) < -d'(r)$, in light of the above and \eqref{10}, it suffices to show:
\begin{equation}
\frac{\mu g\left(\frac{1}{m_2} \log(1 + m_2 b)\right) \epsilon^{N_-}}{\Lambda R^{N_+ - 1} N_-} > \frac{lmb}{R - \epsilon}.
\end{equation}

Set $\epsilon_0 = \frac{N_- R}{N_- + 1}$ where the function $\frac{1}{(R - \epsilon) \epsilon^{N_-}}$ is minimized. From condition $(C2)$, we have:
\[
\mu > \frac{b}{g\left(\frac{1}{m_2} \log(1 + m_2 b)\right)} \frac{\Lambda N_- R^{N_+ - 1}}{(R - \epsilon_0) \epsilon_0^{N_-}}.
\]
Thus, selecting $l, m > 1$ such that:
\[
\mu > \frac{lmb \Lambda N_- R^{N_+ - 1}}{g\left(\frac{1}{m_2} \log(1 + m_2 b)\right) (R - \epsilon_0) \epsilon_0^{N_-}}
\]
ensures that $\tilde{\psi}_2 = \frac{1}{m_2} \log(1 + m_2 \psi_2)$ is a strict subsolution for $\mu > \mu_*$.

Since $\psi_2$ at the origin exceeds the corresponding value of $d$, namely $\psi_2(0) > d(0) = b \rho(0) = b$, we deduce that 
\[
\tilde{\psi}_2(0) = \frac{1}{m_2} \log(1 + m_2 b) > \frac{1}{m_2} \log(1 + m_2 \psi_2(0)).
\]
Given that $b > a$, we have
\[
\|\tilde{\phi}_2\|_{\infty} < \frac{1}{m_1} \log(1 + m_1 a),
\]
which implies
\[
\frac{1}{m_2} \log(1 + m_2 b) > \frac{1}{m_1} \log(1 + m_1 a).
\]
Thus, $\tilde{\psi}_2 \nleq \tilde{\phi}_2$.

According to Theorem \ref{Th4.2}, we can select a sufficiently small subsolution $\tilde{\psi}_1$ and a sufficiently large supersolution $\tilde{\phi}_1$ such that
\[
\tilde{\psi}_1 \leq \tilde{\psi}_2 \leq \tilde{\phi}_1 \quad \text{and} \quad \tilde{\psi}_1 \leq \tilde{\phi}_2 \leq \tilde{\phi}_1.
\]
By Theorem \ref{Th3.7}, these relationships ensure the existence of at least three positive solutions for $\mu$ in the interval $(\mu_*, \mu^*)$, where $\mu^*$ is defined as
\[
\mu^* = \min\left\{\frac{a}{\|e\|_{\infty} (1 + m_1 a) g(a)}, \hat{\mu}\right\}.
\]

For a general bounded domain $\Omega$, let $B_R$ denote the largest inscribed ball in $\Omega$. To demonstrate the existence of three positive solutions, we employ these sub and supersolutions in $B_R$. In $\Omega \setminus B_R$, we extend by zero and verify that at least three positive solutions exist for $\mu \in (\mu_*, \mu^*)$ in the general bounded domain.
\end{proof}


\section{Example}
\label{sec:example}

We illustrate the application of Theorems \ref{Th4.1}-\ref{Th4.3} with a specific example. Consider the boundary value problem:
\begin{equation}\label{exp}
\begin{aligned}
-\mathcal{M}_{\lambda,\Lambda}^+(D^2u) - \Gamma|Du|^2 &= \mu\left(\frac{2u^{1-\alpha}}{1 + u^{1-\alpha}} + e^{\frac{\tau u}{\tau + u}} - 1\right) & &\text{in } \Omega,\\
u &= 0 & &\text{on } \partial\Omega,
\end{aligned}
\end{equation}
where $0 < \alpha < 1$ and $\tau > 0$. Define the function $g(u) = \frac{2u^{1-\alpha}}{1 + u^{1-\alpha}} + e^{\frac{\tau u}{\tau + u}} - 1$.

\begin{itemize}
    \item The function $g(u)$ satisfies $g(0) = 0$ and is an increasing function of $u$. Additionally, we observe that $\lim_{s \to \infty} g(s) = e^{\tau} + 1$ and $\lim_{s \to 0} \frac{g(s)}{s} = \infty$. These properties confirm that the problem \eqref{exp} meets conditions $C0-C2$.
    \item For a sufficiently large $\tau$, choosing $a = 1$ and $b = \tau$, the critical value $\mu_*$ approaches zero as $\mu_* = A \frac{\tau}{g\left(\frac{1}{m_2} \log(1 + m_2 \tau)\right)}$. This leads to $\mu_* < \frac{1}{\|e\|_{\infty}(1 + m_1) e^{\frac{\tau}{\tau + 1}}}$ and $\mu_* < \frac{1}{\|e\|_{\infty}(1 + m_1) (e^{\tau} + 1)} = \hat{\mu}$.
\end{itemize}

From Theorem \ref{Th4.2}, it follows that for $\mu$ values in the range $0 \leq \mu \leq \hat{\mu}$, problem \eqref{exp} admits at least one positive solution. Expanding upon this, Theorem \ref{Th4.3} ensures that for $\mu$ in the interval $(\mu_*, \mu^*)$, where
\[
\mu^* = \min\left\{\frac{1}{\|e\|_{\infty}(1 + m_1) e^{\frac{\tau}{\tau + 1}}}, \frac{1}{\|e\|_{\infty}(1 + m_1) (e^{\tau} + 1)}\right\},
\]
there exist at least three positive solutions. Additionally, Theorem \ref{nonexist} indicates that for sufficiently large $\mu$, no positive solutions exist for \eqref{exp}.

\bibliography{Multiplicity_gradient.bib}
\bibliographystyle{abbrv}
\newpage
\end{document}